\documentclass{amsart}

\newtheorem{theorem}{Theorem}[section]
\newtheorem{lemma}[theorem]{Lemma}

\theoremstyle{definition}
\newtheorem{definition}[theorem]{Definition}
\newtheorem{example}[theorem]{Example}

\theoremstyle{remark}
\newtheorem{remark}[theorem]{Remark}

\numberwithin{equation}{section}

\begin{document}

\title{A new understanding of $\zeta(k)$}

\author{Chenfeng He}

\address{school of Mathematical Sciences, Zhejiang University, Hangzhou
, China 310027}

\begin{abstract}
In this paper, by introducing a new operation in the vector space of Laurent series, the author derived  explicit series for the values of $\zeta$-funtion at positive integers, where $\zeta$ denotes the Riemann zeta function. The values of $\zeta(k),\  k>1$ are largely connected with Bernoulli numbers and binomial numbers. The method in this paper seems new, and the resluts are about divergent series. Using Borel summation for these divergent series one can connect $\zeta$ function, Bernoulli numbers, and  most series representations of Riemann zeta function.
 \end{abstract}
\maketitle

\section{introduction}
In \cite{he2018new} we introduced an operation called convolution, denoted $\star$ in the vector space of Taylor series. If an analytic function $f(z)=\sum_{n=0}^{\infty}a_n\frac{z^{n}}{n!}$, then we denote it by $\boldsymbol{f}=(a_0,a_1,a_2,\cdots)$. The operation is \[  \boldsymbol{f}\star\boldsymbol{g}=(a_0b_0,a_1b_0+a_0b_1,\cdots,\sum_{k=0}^{n}\binom{n}{k}a_kb_{n-k},\cdots),\] where $\boldsymbol{g}=(b_0,b_1,b_2,\cdots)$.
This is actually using  $\{\frac{z^{n}}{n!} \mid n=0,1,2,3,...\}$ as a basis to represent analytic functions and we use the operation to deduce $\zeta(-k)=-\frac{B_{k+1}}{k+1},$ where $k\in \mathbb{N}$. 

Now we generalize this operation to Laurent series to deduce  an ``absurd" series for values of $\zeta$-function at positive integers. That is \[(n-1)(\zeta(n)-1)=\sum_{k=0}^{\infty}\binom{n-2+k}{n-2}B_k.\]
Then we turn ``absurd" to rigorous by Borel summation and give a new understanding of $\zeta$-function and Bernoulli numbers.

\subsection{Laurent series}
The Laurent series is of the form $\sum_{n=-\infty}^{\infty}a_n(z-c)^n$, and it can  represent a complex function $f(z)$ by a power series which includes terms of negative degree. Assume that $f$ is holomorphic in an annulus $A(c,R_1,R_2)$, then $f$ is the sum of a uniquely determined Laurent series: \[f(z)=\sum_{n=-\infty}^{\infty}a_n(z-c)^n,\] where  $a_n$ is defined by\[a_n=\frac{1}{2\pi i}\int_{|z-c|=r'}\frac{f(z)}{(z-c)^{n+1}}dz.\] The circle $|z-c|=r$ lies in the annulus $\boldmath{A}$ (see\cite{ahlfors1953complex}), where $A$ is $\{z: R_1<|z-c|<R_2\}.$ 
 In \cite{he2018new} we use $\{\frac{z^{n}}{n!}\mid n\in \mathbb{N}\}$ as a basis of vector space which consists of analytic functions. Now we generalize this idea. 
\begin{definition}

	Consider all the complex functions holomorphic in some annulus $A(0,r,R),$ which make a vector space over $\mathbb{C}$. They can be uniquely represented as a Laurent series. We use \[\{\frac{(-1)^{n-1}(n-1)!}{z^{n}}, 1, \frac{z^{n}}{n!}\mid n=1,2,3,...\}\] as a basis of this vector space. If \[f(z)=\sum_{n=1}^{\infty}a_{-n}\frac{(-1)^{n-1}(n-1)!}{z^{n}}+\sum_{n=0}^{\infty}a_n\frac{z^{n}}{n!}\] in the domain $A(0,r_1,R_1)$, then we denote $f(z)$ by $\boldsymbol{f}=(\cdots,a_{-2},a_{-1}\mid a_0,a_1,a_2,\cdots).$ Assuming another function $\boldsymbol{g}=(\cdots,b_{-2},b_{-1}\mid b_0,b_1,b_2,\cdots)$ in the annulus $A(0,r_2,R_2)$ and the intersection of their domains is not empty, we can define addition: \[f(z)+g(z)=\boldsymbol{f}+\boldsymbol{g}=(\cdots,a_{-2}+b_{-2},a_{-1}+b_{-1}\mid a_0+b_0,a_1+b_1,a_2+b_2,\cdots);\] and scalar multiplication: \[\lambda f(z)=\lambda \boldsymbol{f}=(\cdots,\lambda a_{-2},\lambda a_{-1}\mid \lambda a_0,\lambda a_1,\lambda a_2,\cdots), \] where $\lambda \in \mathbb{C}$.  We use a symbol  $\boldsymbol{\cdot}$ to indicate function multiplication, i.e. $f(z)g(z)=\boldsymbol{f\cdot g}$.
	
	If $\boldsymbol{f}=(\cdots,a_{-2},a_{-1}\mid a_0,a_1,a_2,\cdots),$ we use $(\boldsymbol{f})_{\pm n}=a_{\pm n}$ to denote the $\pm n$-th component of $\boldsymbol{f}$. 
	
	Now we define for $n\geq 0$ \begin{equation} (\boldsymbol{f\cdot g})_{n}=\sum_{k=-\infty}^{-1}\frac{(-1)^{k}}{k}\frac{1}{\binom{n-k}{n}}a_kb_{n-k}+\sum_{k=0}^{n}\binom{n}{k}a_kb_{n-k}+\sum_{k=n+1}^{\infty}\frac{(-1)^{n-k}}{n-k}\frac{1}{\binom{k}{n}}a_kb_{n-k},\end{equation} and for $n>0$ \begin{equation}
	(\boldsymbol{f\cdot g})_{-n}=\sum_{k=-\infty}^{-n}(-1)^{n+k}\binom{-k-1}{n-1}a_kb_{-n-k}+\sum_{k=-n+1}^{-1}\frac{1}{k}\frac{1}{\binom{n-1}{-k}}a_kb_{-n-k}+\sum_{k=0}^{\infty}(-1)^{k}\binom{n+k-1}{n-1}a_kb_{-n-k}.\end{equation} \end{definition} 
This is obtained by comparing the $n$-th component of $f(z)g(z)$. For simplicity, we need to generalize the binomial number.
\begin{definition} For $p, q\geq0$, we define 
\begin{equation}
\left< \begin{matrix}
p\\
q
\end{matrix}\right> =	
\begin{cases} 
\binom{p}{q},  &\mbox{if }p\geq q\\
\frac{(-1)^{p-q}}{p-q}\frac{1}{\binom{q}{p}} & \mbox{if} p<q
\end{cases}.
\end{equation}
\end{definition}
Therefore (1.1) becomes \begin{equation} (\boldsymbol{f\cdot g})_{n}=\sum_{k=-\infty}^{-1}\left< \begin{matrix}
n\\
n-k
\end{matrix}\right>a_kb_{n-k}+\sum_{k=0}^{n}\left< \begin{matrix}
n\\
k
\end{matrix}\right>a_kb_{n-k}+\sum_{k=n+1}^{\infty}\left< \begin{matrix}
n\\
k
\end{matrix}\right>a_kb_{n-k}.\end{equation} 
\begin{lemma}By  definition, $\boldsymbol{\cdot}$ has the following rules:
\begin{enumerate}
\item \text{Commutative}: 
$\boldsymbol{f\cdot g}=\boldsymbol{g\cdot f}$
\item \text{Associative} : $
\boldsymbol{(f\cdot g)\cdot h}=\boldsymbol{f\cdot (g\cdot h)}$

\item \text{Distributive}: $
\boldsymbol{f\cdot}\boldsymbol{(g}+\boldsymbol{h)}=\boldsymbol{f\cdot g}+\boldsymbol{f\cdot h}$

\item $\boldsymbol{f\cdot}\lambda\boldsymbol{g}=\lambda\boldsymbol{f\cdot g}=\lambda\boldsymbol{(f\cdot g)}$
\end{enumerate}
\end{lemma}
Because these operations are corresponding to those of  functions, the rules are obvious.
\begin{example}
Let\[e^{z}=\sum_{n=0}^{\infty}\frac{z^{n}}{n!}=(\cdots,0,0,0\mid 1,1,1,\cdots)\] and \[e^{\frac{1}{z}}=\sum_{n=0}^{\infty}\frac{z^{-n}}{n!}=(\cdots, \frac{-1}{3!\times 4!},\frac{1}{2!\times 3!},\frac{-1}{1!\times2!},1\mid 1,0,0\cdots),\]
obviously the intersection of their domains is not empty, so \[e^{z}+e^{\frac{1}{z}}=(\cdots,\frac{-1}{3!\times 4!},\frac{1}{2!\times 3!},\frac{-1}{1!\times2!},1\mid 2,1,1,\cdots),\] and \[e^{z}e^{\frac{1}{z}}=\sum_{n=1}^{\infty}\sum_{k=0}^{\infty}\frac{1}{k!(n+k)!}z^{-n}+\sum_{n=0}^{\infty}(1+\sum_{k=n+1}^{\infty}\frac{n!}{(k-n)!k!})z^{n}.\]\end{example} If we define the identity to be $\boldsymbol{id}=(\cdots,0,0,0\mid 1,0,0,\cdots)$ and $\boldsymbol{f\cdot g}=\boldsymbol{id}$, we say $ \boldsymbol{g}$ is the inverse of $\boldsymbol{f}$. 

For example, consider $e^{-z}=(\cdots,0,0,0\mid 1,-1,1,-1,\cdots)$, we have \[ e^{z}e^{-z}=(\cdots,0,0,0\mid 1,1,1,\cdots)\boldsymbol{\cdot}(\cdots,0,0,0\mid 1,-1,1,-1\cdots)=\boldsymbol{id},\] hence the inverse of $e^{z}$ is $e^{-z}$.

\subsection{Vector multiplication} Define for $j\in \mathbb{C},\ j\neq0$ 
\[\boldsymbol{j}=(\cdots,0,0\mid j^0,j^1,j^2,\cdots), \ \boldsymbol{j^{-1}}=(\cdots,0,0,0\mid(-j)^0,(-j)^1,(-j)^2,\cdots).\]  One can see $\boldsymbol{j}=e^{jz},\  \boldsymbol{j^{-1}}=e^{-jz}$, so we get the inverse of $\boldsymbol{j}$ is $\boldsymbol{j^{-1}}$. 

 We also define vector multiplication to be \[\boldsymbol{f}\boldsymbol{g}=(\cdots,a_{-2}b_{-2},a_{-1}b_{-1}\mid a_0b_0,a_1b_1,a_2b_2,\cdots).\] 
Therefore we can write $\boldsymbol{j^{-1}}=\boldsymbol{-1}\boldsymbol{j}.$

\subsection{Bernoulli numbers}
Recall that the generating function of Bernoulli numbers is
 \begin{equation}
\frac{z}{e^z-1}=\sum_{k=0}^{\infty}B_k\frac{z^k}{k!},\end{equation} and we denote it by $\boldsymbol{B}=(\cdots,0,0,0\mid B_0,B_1,B_2,\cdots)$. The inverse of $\boldsymbol{B}$ is \begin{equation}
\frac{e^z-1}{z}=\sum_{k=0}^{\infty}\frac{1}{k+1}\frac{z^k}{k!},
\end{equation}
 denoted $\boldsymbol{H}=(\cdots,0,0,0\mid 1,\frac{1}{2},\frac{1}{3},\cdots)$. 
 We write \begin{equation}
 \boldsymbol{B\cdot H}=\boldsymbol{id}.
\end{equation} 
We also note that \[\frac{-z}{e^{-z}-1}=\frac{ze^{z}}{e^{z}-1},\]which means \begin{equation}
\boldsymbol{-1}\boldsymbol{B}=\boldsymbol{B\cdot 1}.\end{equation}
Moreover \[\frac{e^{-z}-1}{-z}=\frac{e^{z}-1}{ze^{z}},\] hence we have \begin{equation}
\boldsymbol{-1}\boldsymbol{H}=\boldsymbol{H\cdot -1}.
\end{equation}
We also have \begin{equation}
\boldsymbol{-1B\cdot -1H}=\boldsymbol{id}.
\end{equation}
We will use these identities in the next section.
\section{$\zeta$-function}
We need to be clear that the coefficents of $\frac{(-1)^{n-1}(n-1)!}{z^{n}}$ are $0$, the Laurent series becomes Taylor series, and the $\boldsymbol{\cdot}$ becomes $\star$. Then all the relations in \cite{he2018new} apply.   
\subsection{}Riemann zeta function is defined as \[\zeta(x)=\sum_{n=1}^{\infty}\frac{1}{n^s}\] for $Re(s)>1$, and extends to an analytic function for all $s\in \mathbb{C}$, except for $s=1$, where it has a simple pole. The way of extension can be easily seen in \cite{apostol2013introduction} and \cite{murty2000simple} as is without using functional equation. The conclution is \begin{equation}
\zeta(s)=1+\frac{1}{s-1}-\sum_{r=1}^{m}\frac{s(s+1)\cdots(s+r-1)}{(r+1)!}(\zeta(s+r)-1)-\frac{s(s+1)\cdots(s+m)}{(m+1)!}\sum_{n=1}^{\infty}\int_{0}^{1}\frac{u^{m+1}du}{(u+n)^{s+m+1}} \end{equation}
and the sum on the right hand converges for $Re(s)>-m$. In \cite{he2018new} we use this to deduce 
 	\begin{equation} \delta_{m,0}-\frac{(-1)^{m}}{m+1}=\sum_{i=0}^{m}(1-\delta_{m,0})(-1)^{i}\binom{m}{i}\zeta(i-m).  \end{equation}
 This  is equivalent to \begin{equation}
 \boldsymbol{id}-(\boldsymbol{-1H})=\boldsymbol{-1H\cdot}(\cdots,0,0\mid 0,-1,0,0,\cdots)\boldsymbol{\cdot \zeta(-s)},
 \end{equation}
  where $\boldsymbol{ \zeta(-s)}=(\cdots,0,0,0\mid\zeta(0),\zeta(1),\zeta(2),\cdots).$ 
For $n\geq0$ we get the well-known relaton \begin{equation}
\zeta(-n)=-\frac{B_{n+1}}{n+1}.\end{equation}
\subsection{}What would happen when we put positive integers into (2.1)?
 \cite{M1994remark} mentioned a beautiful formula which is   \begin{equation}
(s-1)(\zeta(s)-1)-1=-\sum_{r=1}^{\infty}\frac{(s-1)s\cdots(s+r-1)}{(r+1)!}(\zeta(s+r)-1).
\end{equation}
This also can be deduced from (2.1) by multiplying $(s-1)$ to both sides, and letting $m\to\infty$.
Let's take a look at (2.5), after moving terms we get \begin{equation}
\sum_{r=0}^{\infty}\frac{(s-1)s\cdots(s+r-1)}{(r+1)!}(\zeta(s+r)-1)=1.
\end{equation}
Now we can use our operation.

Define $\boldsymbol{\zeta'(s)}$ to be $(\cdots,3(\zeta(4)-1),2(\zeta(3)-1),\zeta(2)-1\mid 0,0,\cdots)$, and we change (2.6) to
\begin{equation}
\sum_{r=0}^{\infty}\frac{(-1)^{r}}{r+1}(-1)^{r}\binom{s+r-2}{r}(s+r-1)(\zeta(s+r)-1)=1.
\end{equation}
This can be regraded as the negative component of \begin{equation}
(\cdots,0,0,0\mid 1,-\frac{1}{2},\frac{1}{3},\cdots)\boldsymbol{\cdot}(\cdots,3(\zeta(4)-1),2(\zeta(3)-1),\zeta(2)-1\mid 0,0,\cdots).
\end{equation}
 We know $(\cdots,0,0,0\mid 1,-\frac{1}{2},\frac{1}{3},\cdots)=\boldsymbol{-1H}$, we can denote (2.7) as 
 \begin{equation}
 \boldsymbol{-1H\cdot\zeta'(s)}=(\cdots,1,1,1\mid a_0,a_1,a_2,\cdots).
 \end{equation}
  We multiply $\boldsymbol{-B}$ on both side of (2.9), and from (1.10) we get 
 \begin{equation}
 \boldsymbol{id\cdot\zeta'(s)}=\boldsymbol{-B\cdot}(\cdots,1,1,1\mid a_0,a_1,a_2,\cdots).
 \end{equation}
 Making use of (1.2) to compare the negative components of (2.10), we get 
 \begin{theorem} For positive integers $n\geq 2$
\begin{equation}
(n-1)(\zeta(n)-1)=\sum_{k=0}^{\infty}\binom{n-2+k}{n-2}B_k.
\end{equation} 
 \end{theorem}
The reader may notice that (2.11) doesn't make sense, because the infinite summation on the right hand obviously diverge. The term $|\binom{n-2+k}{n-2}B_k|$ does not go to $0$ as $k\rightarrow \infty$. We will explain theorem 2.1 in the next section.
 
\section{divergent series}
Maybe the most famous divergent series was Grandi's series, which is
\begin{equation}
1-1+1-1+\cdots=\sum_{n=0}^{\infty}(-1)^{n}.
\end{equation}
This series was reported by Guido Grandi in 1703(see the history in \cite{bagni2009inf}). By inserting parentheses into $1-1+1-1+\cdots$, we produce different results: either
$(1-1)+(1-1)+\cdots=0$ or $1+(-1+1)+(-1+1)+\cdots=1.$

In the 1700's, many mathematicians wanted to find a value for this series and they didn't think (3.1) summed to either 0 or 1. Actually most of them thought the true value is $\frac{1}{2}$. There are many explanations about this value. For example, Daniel Bernoulli thought that since half of the partial sums of (3.1) are $+1$ and half of them are $0$, the correct value of the series would be $\frac{1}{2}\times1+\frac{1}{2}\times0=\frac{1}{2}$. The usual reason for the value lies in putting $x=-1$ in the following geometric series
\begin{equation}\frac{1}{1-x}=1+x+x^{2}+x^{3}+\cdots .\end{equation}
But we know (3.2) converges only if $|x|<1$. All of these do not make a rigorous way to explain the divergent series.

Divergent series were widely used by Leonhard Euler, but often led to confusing and contradictory results. His idea that any divergent series should have a natural sum had been hiding in the sea of mathematics, since Cauchy gave a rigorous definition of the sum of a convergent series. They reappeared in 1886 with Henri Poincaré's work on asymptotic series. In 1890, Ernesto Cesàro gave a rigorous definition of the sum of some divergent series(like Grandi's series), and defined Cesàro summation.

In modern mathematics we already have theorems on methods for summing divergent series. We will use Borel summation to explain theorem 2.1.

\subsection{Borel summation}

Let $A(z)$ denote formal power series 
\[A(z)=\sum_{n=0}^{\infty}a_nz^{n},\] 
and define Borel transform of $A(z)$ to be 
\[BA(t)=\sum_{n=0}^{\infty}\frac{a_n}{n!}t^{n}.\]
\begin{definition}
	Suppose that the Borel transform converges for all positive real numbers to a function  that the following integral is well defined(as an improper integral), the Borel sum of $A$ is given by
	\[\int_{0}^{\infty}e^{-t}BA(zt){\rm d}t.\]
If the integral converges at $z\in \mathbb{C}$ to some $a(z)$, we say that the Borel sum of $A(z)$ converges at $z$, and write \[\sum_{n=0}^{\infty}a_nz^{n}=a(z)(\rm B).\]
\end{definition}

Let's look at some examples.
\begin{example} Let
 $A_1(z)=1+z+z^{2}+z^{3}+\cdots$, then \[BA_1(z)=1+z+\frac{z^{2}}{2!}+\frac{z^3}{3!}+\cdots=e^{z}.\]
The Borel sum of $A_1(z)$  is \[\int_{0}^{\infty}e^{-t}e^{zt}{\rm d}t=\int_{0}^{\infty}e^{(z-1)t}{\rm d}t=\frac{1}{1-z}, \  Re(z)<1.\]
Putting $z=-1,$ we get \[\int_{0}^{\infty}e^{-2t}{\rm d}t=\frac{1}{2}.\]
So the Grandi's series \begin{equation}
1-1+1-1+\cdots=\frac{1}{2}({\rm B}).\end{equation} \\ 
Putting $z=-2$ we have \[\int_{0}^{\infty}e^{-3t}{\rm d}t=\frac{1}{3},\]
that is \begin{equation}
1-2+4-8+\cdots=\frac{1}{3}({\rm B}).\end{equation} \\ 
Putting $z=-\frac{1}{2},$ we have \[\int_{0}^{\infty}e^{-\frac{3}{2}t}{\rm d}t=\frac{2}{3}=1-\frac{1}{2}+(\frac{1}{2})^{2}-(\frac{1}{2})^{3}+\cdots.\]
Putting $z=\frac{1}{2}$, we get \[\int_{0}^{\infty}e^{-\frac{1}{2}t}{\rm d}t=2=1+\frac{1}{2}+(\frac{1}{2})^{2}+(\frac{1}{2})^{3}+\cdots.\]
\end{example}
Noticing that whenever $A(z)$ converges in the standard sense, the Borel sum  converges to the same value, i.e.
\begin{equation}\sum_{n=0}^{\infty}a_nz^n=A(z)<\infty \  \Longrightarrow \sum_{n=0}^{\infty}=a_nz^n=A(z)({\rm B}). \end{equation}
The above property is called regularity of Borel summation method. 	It can be seen by a change in the order of integration, which is valid due to absolute convergence: if $A(z)$	is convergent at $z$, then 
\begin{equation}
A(z)=\sum_{n=0}^{\infty}a_nz^n=\sum_{n=0}^{\infty}a_n(\int_{0}^{\infty}e^{-t}t^{n}{\rm d}t)\frac{z^n}{n!}=\int_{0}^{\infty}e^{-t}\sum_{n=0}^{\infty}a_n\frac{(tz)^n}{n!}{\rm d}t,
\end{equation}
where $\int_{0}^{\infty}e^{-t}t^n{\rm d}t=\Gamma(n+1)$.
\begin{example}
	Let $A_2(z)=\sum_{k=0}^{\infty}B_kz^{k}$, where $B_k$ are Bernoulli numbers, then \begin{equation}BA_2(z)=\sum_{k=0}^{\infty}B_k\frac{z^{k}}{k!}=\frac{z}{e^{z}-1}.\end{equation}
	The Borel sum of $A_2(z)$ is \begin{equation}\int_{0}^{\infty}e^{-t}\frac{zt}{e^{zt}-1}{\rm d}t.\end{equation}
	Putting $z=1$, we have \begin{align*}\int_{0}^{\infty}\frac{e^{-t}t}{e^{t}-1}{\rm d}t &=\int_{0}^{\infty}\frac{(e^{-t}-1)t}{e^{t}-1}{\rm d}t+\int_{0}^{\infty}\frac{t}{e^{t}-1}{\rm d}t \\ & =-\int_{0}^{\infty}e^{-t}t{\rm d}t+\Gamma(2)\zeta(2)\\
		&=\Gamma(2)\zeta(2)-\Gamma(2),\end{align*}
then we get \begin{equation}
B_0+B_1+B_2+B_3+\cdots=\zeta(2)-1(\rm B).
\end{equation}
This is exactly a result of Theorem 2.1 when $n=2$.
\end{example}
\begin{example}
Let $A_3(z)=0+B_0z+2B_1z^{2}+3B_2z^{3}+\cdots=\sum_{k=1}^{\infty}kB_{k-1}z^{k}$, then 
\begin{equation}
BA_3(z)=\sum_{k=1}^{\infty}\frac{kB_{k-1}}{k!}z^{k}=z\sum_{k=0}^{\infty}\frac{B_{k}}{k!}z^{k}=\frac{z^{2}}{e^{z}-1}.
\end{equation}
The Borel sum of $A_3(z)$ is 
\begin{equation} \int_{0}^{\infty}\frac{e^{-t}t^{2}z^2}{e^{zt}-1}{\rm d}t.
\end{equation}
Putting $z=1$ in (3.11), we have 
\begin{equation}\int_{0}^{\infty}\frac{e^{-t}t^2}{e^{t}-1}{\rm d}t=\int_{0}^{\infty}\frac{(e^{-t}-1)t^2}{e^{t}-1}{\rm d}t+\int_{0}^{\infty}\frac{t^2}{e^{t}-1}{\rm d}t=\Gamma(3)\zeta(3)-\Gamma(3).
	\end{equation}
This is also a reslut of Theorem 2.1, i.e. \begin{equation}
B_0+2B_1+3B_2+\cdots=\sum_{k=0}^{\infty}(k+1)B_k=2\zeta(3)-2(\rm B).
\end{equation}
\end{example}
The reader may notice that (3.13) doesn't have $0$ as its first term, but $A_3(z)$ does. This is because $0+\sum_{k=0}^{\infty}(k+1)B_k=\sum_{k=0}^{\infty}(k+1)B_k$. Do not take this reason for granted, because it's not like convergent series, in which $0$ equals nothing. In divergent series, $0$ can change a lot, hence we need a rigorous reason for this.

 Actually Borel proved(see\cite{borel1975lectures}) that in a series that is absolutly summable, transposing a finite number of terms  or replacing a certain number of consecutive terms by their sum, or replacing a term by the sum of several others won't change either the summability or the sum of the series. 
 
 If two series are summable by Borel summation method, so does their linear combination, i.e. For $\alpha ,\ \beta \in\mathbb{C},$
 \begin{equation} 
 	A=a_0+a_1+a_2+\cdots=a(B),\ V=v_0+v_1+v_2+\cdots=b(B) \Longrightarrow
 	\alpha A+\beta V=\alpha a+\beta b({\rm B}).
 \end{equation}
\subsection{}Let $A_4(z)=\sum_{k=0}^{\infty}(k+1)B_kz^{k}=\sum_{k=0}^{\infty}B_kz^{k}+\sum_{k=0}^{\infty}kB_kz^{k}$, then \[BA_4(z)=\sum_{k=0}^{\infty}\frac{B_k}{k!}z^{k}+\sum_{k=0}^{\infty}\frac{kB_k}{k!}z^{k}=\frac{z}{e^{z}-1}+z(\frac{z}{e^{z}-1})'.\]
The Borel sum of $A(z)$ is 
\begin{align*} &\int_{0}^{\infty}\frac{e^{-t}zt}{e^{zt}-1}{\rm d}t+\int_{0}^{\infty}\frac{e^{-t}zt(e^{zt}-1-zte^{tz})}{(e^{tz}-1)^{2}}{\rm d}t\\ &=2\int_{0}^{\infty}\frac{e^{-t}zt}{e^{zt}-1}{\rm d}t-\int_{0}^{\infty}\frac{e^{-t}(zt)^{2}e^{tz}}{(e^{tz}-1)^{2}}{\rm d}t \\ &=2\int_{0}^{\infty}\frac{e^{-t}zt}{e^{zt}-1}{\rm d}t-\int_{0}^{\infty}\frac{(e^{-t}-1+1)(zt)^{2}e^{zt}}{(e^{zt}-1)^{2}}{\rm d}t.
\end{align*}
Putting $z=1$ in the above formula, we have 
\begin{align}&2\int_{0}^{\infty}\frac{e^{-t}t}{e^{t}-1}{\rm d}t+\int_{0}^{\infty}\frac{t^{2}}{e^{t}-1}{\rm d}t-\int_{0}^{\infty}\frac{e^{t}t^{2}}{(e^{t}-1)^{2}}{\rm d}t\\ &=2\zeta(2)-2+\Gamma(3)\zeta(3)-\int_{0}^{\infty}\frac{e^{t}t^{2}}{(e^{t}-1)^{2}}{\rm d}t.
\end{align}
From \cite{Ragib2012evaluating}, the last term of (3.16) is $2\zeta(2)$. Therefore (3.16) becomes $\Gamma(3)\zeta(3)-2$, this is to say 
\begin{equation}
\sum_{k=0}^{\infty}(k+1)B_k=2\zeta(3)-2(\rm B). 
\end{equation}  Comparing with example 3.4, we have\[0+\sum_{k=0}^{\infty}(k+1)B_k=\sum_{k=0}^{\infty}(k+1)B_k.\]
\subsection{proof of Theorem 2.1}
First, we multiply $(n-2)!$ on both sides of (2.11) and get
\begin{equation}
\Gamma(n)(\zeta(n)-1)=\sum_{k=0}^{\infty}(n-2+k)_{n-2}B_k(\rm B),
\end{equation}
where $(n-2+k)_{n-2}$ is falling factorial.

Second, let $D_n(z)=\sum_{i=1}^{n-2}0+\sum_{k=0}^{\infty}(n-2+k)_{n-2}B_kz^{n-2+k}$ where $n>2$ and $D_2(z)=\sum_{k=0}^{\infty}B_kz^{k}$. Then for positive integers $n\geq2$, we get \begin{equation}
BD_n(z)=\sum_{k=0}^{\infty}\frac{(n-2+k)_{n-2}B_k}{(n-2+k)!}z^{n-2+k}
=\sum_{k=0}^{\infty}\frac{B_k}{k!}z^{n-2+k}=\frac{z^{n-1}}{e^{z}-1}.
\end{equation}

Third, the Borel sum of $D_n(z)$ is 
\begin{align} \int_{0}^{\infty}\frac{e^{-t}(zt)^{n-1}}{e^{zt}-1}{\rm d}t=\int_{0}^{\infty}\frac{(e^{-t}-1)(zt)^{n-1}}{e^{zt}-1}{\rm d}t+\int_{0}^{\infty}\frac{(zt)^{n-1}}{e^{zt}-1}{\rm d}t.\end{align}

Last, putting $z=1$ in (3.20), we get 
\begin{equation}
\int_{0}^{\infty}\frac{(e^{-t}-1)t^{n-1}}{e^{t}-1}{\rm d}t+\int_{0}^{\infty}\frac{t^{n-1}}{e^{t}-1}{\rm d}t=\Gamma(n)\zeta(n)-\Gamma(n).\end{equation}
\qed

 If we put $z=-1$ in (3.20), we get \begin{equation}
\int_{0}^{\infty}\frac{e^{-t}(-t)^{n-1}}{e^{-t}-1}{\rm d}t=(-1)^{n}\int_{0}^{\infty}\frac{t^{n-1}}{e^{t}-1}{\rm d}t=(-1)^{n}\Gamma(n)\zeta(n).
\end{equation}	
This is to say 
\begin{equation}
\sum_{k=0}^{\infty}(-1)^{n-2}(n-2+k)_{n-2}(-1)^{k}B_k=(-1)^{n}\Gamma(n)\zeta(n)(\rm B).
\end{equation}
Therefore we have 
\begin{theorem} For positive integers $n>1$,
	 \begin{equation}
	\sum_{k=0}^{\infty}(n-2+k)_{n-2}(-1)^{k}B_k=\Gamma(n)\zeta(n)(\rm B).\end{equation}
\end{theorem}
\subsection{Algebra structure}
Actually Borel also proved there is a multiplication among absolutly summable series, i.e. if 
\[W=w_0+w_1+w_2+\cdots=w({\rm{B}}), \ V=v_0+v_1+v_2+\cdots=v({\rm{B}}),\]
then \begin{align}W\cdot V=\sum_{n=0}^{\infty}c_n=wv(\rm B),\end{align} where $c_n=\sum_{k=0}^{n}w_kv_{n-k}$.

If we define any absolutely summable series $A(z)=\sum_{n=0}^{\infty}a_nz^{n}$ as ordered sequence $(a_0,a_1,a_2,a_3,\cdots)$, then they form an algebra over $\mathbb{C}$. If a series $A(z)$ is summable at $z_0$ and has a value $a$, then we write \begin{equation}
A(z_0)=(a_0,a_1,a_2,a_3,\cdots)_{z_0}=a({\rm B}).
\end{equation}
Let's look at some examples that a series isn't summable at some $z_0$.
\begin{example}
Recalling that $A_1(z)=(1,1,1,1,\cdots)$ in example 3.2, we know that $(1,1,1,\cdots)_{-1}=\frac{1}{2}({\rm B})$. What about $(1,1,1,\cdots)_1$? The Borel sum of $A_1(1)$ is \[\int_{0}^{\infty}e^{-t}e^{t}{\rm d}t=\infty,\] therefore $1+1+1+\cdots$ is not summable by Borel's method.
\end{example}
\begin{example}
Let $N(z)=\sum_{n=0}^{\infty}nz^{n}$, then \begin{equation}
BN(z)=\sum_{n=0}^{\infty}\frac{n}{n!}z^{n}=ze^{z}.
\end{equation}
The Borel sum of $BN(-1)$ is \begin{align}
	\int_{0}^{\infty}e^{-t}(-1)te^{-t}{\rm d}t=-\int_{0}^{\infty}te^{-2t}{\rm d}t=-\frac{1}{4}, \end{align}
which means $0-1+2-3+4-5+\cdots=-1+2-3+4-5+\cdots=-\frac{1}{4}({\rm B}),$
or $1-2+3-4+\cdots=\frac{1}{4}({\rm B}).$ We also note that \begin{equation}(1-1+1-1+\cdots)\cdot(1-1+1-1+\cdots)=1-2+3-4+\cdots=\frac{1}{4}({\rm B}).\end{equation}
The Borel sum of $BN(1)$ is \begin{equation}
\int_{0}^{\infty}e^{-t}te^{t}{\rm d}t=\int_{0}^{\infty}t{\rm d}t=\infty,
\end{equation} therefore $1+2+3+4+\cdots$ is not summable. These results are not like $1+1+1+\cdots=-\frac{1}{2}=\zeta(0)$ and $1+2+3+4+\cdots=-\frac{1}{12}=\zeta(-1)$, which should not be mixed up.
\end{example}
There's one thing we need to be clear: if one side is a divergent series and the other is a finite value, when we move a term from one side to the other, we need to follow a rule which the author would call ``Borel sum rule''. Let's look at an example to illustrate this.
\begin{example} We already had $B_0+B_1+B_2+\cdots=\zeta(2)-1$(B). If someone wants to move $-1$ to the left side, we need to resume $-1$ to the original appearence first. We can think of $-1$ as $-\int_{0}^{\infty}e^{-t}{\rm d}t$, then we move it to the left and get \begin{equation}
2+B_1+B_2+B_3+\cdots=\zeta(2)({\rm B}).
	\end{equation} If we think of $-1$ as $-\int_{0}^{\infty}e^{-t}t{\rm d}t$, then we move it and get \begin{equation}
B_0+(B_1+1)+B_2+B_3+\cdots=B_0-B_1+B_2+B_3+\cdots=\zeta(2){(\rm B)}.
	\end{equation}
If we treat $\Gamma(n)$ as $(n-1)!\int_{0}^{\infty}e^{-t}t{\rm d}t$ then from (3.18) we get 
\begin{equation}
(n-2)!B_0+(n-1)!(B_1+1)+\sum_{k=2}^{\infty}(n-2+k)_{n-2}B_k=\Gamma(n)\zeta(n)({\rm B}).
\end{equation}
This is exactly Theorem 3.5 because $B_{2k+1}=0,\ k>0$. 
\end{example}
\subsection{} There are so many identities related to Bernoulli numbers and Riemman zeta function that sometimes we wonder why these happen.
For instance, from \cite{gessel2005miki} and \cite{agoh2014miki} we have the following equations, for $ n\geq4$
 \begin{align} &\sum_{k=2}^{n-2}\binom{n}{k}B_kB_{n-k}=-(n+1)B_n, \\
 		&(n+2)\sum_{k=2}^{n-2}B_kB_{n-k}-2\sum_{k=2}^{n-2}\binom{n+2}{k}B_kB_{n-k}=n(n+1)B_n, \\  
	&\sum_{k=2}^{n-2}\beta_k\beta_{n-k}-\sum_{k=2}^{n-2}\binom{n}{k}\beta_k\beta_{n-k}=2H_{n}\beta_n,
\end{align}
where $\beta_n=\frac{B_{n}}{n}, \ H_{n}=1+\frac{1}{2}+\cdots+\frac{1}{n} $.
They all can be explained in our divergent series. First we need to add $B_0,B_n,B_1,B_{n-1}$ to these equations.
For the first equation, we change it to \begin{equation}
\sum_{k=0}^{n}\binom{n}{k}B_kB_{n-k}=-(n+1)B_n+2B_0B_n+2nB_1B_{n-1}=-(n-1)B_{n}-nB_{n-1},
\end{equation}
where $n\geq1$.
Now we consider the sequences which is defined in 3.4, i.e. 
\begin{align}&S_1(z)=(0,\sum_{k=0}^{1}\binom{1}{k}B_{k}B_{1-k},\cdots,\sum_{k=0}^{n}\binom{n}{k}B_{k}B_{n-k},\cdots)=\sum_{n=1}^{\infty}\sum_{k=0}^{n}\binom{n}{k}B_{k}B_{n-k}z^{n}, \\
	&S_2(z)=(0,0,-B_2,-2B_3,\cdots,-(n-1)B_n,\cdots)=-\sum_{n=1}^{\infty}(n-1)B_{n}z^{n}, \\
	&S_3(z)=(0,-B_0,-2B_1,\cdots,-nB_{n-1},\cdots)=-\sum_{n=1}^{\infty}nB_{n-1}z^{n}.
\end{align}
As usual we get the Borel transform of (3.38)
\begin{align}BS_1(z)=\sum_{n=1}^{\infty}\sum_{k=0}^{n}\binom{n}{k}B_{k}B_{n-k}\frac{z^{n}}{n!}=(\frac{z}{e^{z}-1})^{2}-1.
\end{align}
Its Borel sum at $z=1$ is \begin{align*}&\int_{0}^{\infty}e^{-t}(\frac{t^{2}}{(e^{t}-1)^{2}}-1){\rm d}t=\int_{0}^{\infty}\frac{(e^{-t}-1+1)t^{2}}{(e^{t}-1)^{2}}{\rm d}t-1\\ 
&=-\int_{0}^{\infty}\frac{e^{-t}t^{2}}{e^{t}-1}{\rm d}t+\int_{0}^{\infty}\frac{(e^{-t}-1+1)t^{2}e^{t}}{(e^{t}-1)^{2}}{\rm d}t-1 \\ 
&=-2\zeta(3)+1-\int_{0}^{\infty}\frac{t^{2}}{e^{t}-1}{\rm d}t+\int_{0}^{\infty}\frac{t^{2}e^{2}}{(e^{t}-1)^{2}}{\rm d}t \\
&=-2\zeta(3)+1-2\zeta(3)+2\zeta(2),
\end{align*}
which means \begin{equation}\sum_{n=1}^{\infty}\sum_{k=0}^{n}\binom{n}{k}B_{k}B_{n-k}=-4\zeta(3)+2\zeta(2)+1({\rm B}).\end{equation}
From Theorem 2.1 we  know that $\sum_{n=1}^{\infty}nB_{n-1}=2\zeta(3)-2$(B), and \[\sum_{n=1}^{\infty}(n-1)B_{n}=\sum_{n=0}^{\infty}(n+1)B_{n}-2\sum_{n=0}^{\infty}B_{n}+1=2\zeta(3)-2\zeta(2)+1({\rm B}).\]
Combining all these, we have \begin{equation}
\sum_{n=1}^{\infty}\sum_{k=0}^{n}\binom{n}{k}B_{k}B_{n-k}=-\sum_{n=1}^{\infty}nB_{n-1}-\sum_{n=1}^{\infty}(n-1)B_{n}.
\end{equation}
From now on we denote $B^{+}_k$ as the second kind of Bernoulli numbers, i.e. $B^{+}_k=B_k$ when $k\neq1$, $B^{+}_1=-B_1.$ After some computation (3.34) becomes \begin{equation}
\sum_{k=0}^{n}\binom{n}{k}B_{k}^{+}B_{n-k}^{+}=-(n-1)B_{n}^{+}+nB_{n-1}^{+},\ n\geq1.
\end{equation}
We can get the similar formula
\begin{equation}
\sum_{n=0}^{\infty}\sum_{k=0}^{n}\binom{n}{k}B_{k}^{+}B_{n-k}^{+}=2\zeta(2)({\rm B}).
\end{equation}

Formula (3.35) was discovered by Matiyasevich\cite{matiyasevich1997identities}. Let's see what this means in divergent series.
First we need to add $B_0^{+},B_n^{+},B_1^{+},B_{n-1}^{+}$ to the formula, and get 
\begin{equation}
(n+2)\sum_{k=0}^{n}B_k^{+}B_{n-k}^{+}+\binom{n+2}{3}B_{n-1}^{+}=2\sum_{k=0}^{n}\binom{n+2}{k}B_k^{+}B_{n-k}^{+},\ n\geq1.
\end{equation}
From Theorem 3.5 we know \[\sum_{k=0}^{\infty}B_{k}^{+}=\zeta(2)({\rm B}),\]and \[\sum_{k=0}^{\infty}(k+1)B_k^{+}=2\zeta(3)({\rm B}).\]
Mutiplying the above we get 
\begin{align}
\sum_{n=0}^{\infty}\sum_{k=0}^{n}B_k^{+}(n-k+1)B_{n-k}^{+}=\sum_{n=0}^{\infty}\sum_{k=0}^{n}(k+1)B_k^{+}B_{n-k}^{+}=2\zeta(2)\zeta(3).
\end{align}
One can see that \[\sum_{k=0}^{n}B_k^{+}(n-k+1)B_{n-k}^{+}+\sum_{k=0}^{n}(k+1)B_k^{+}B_{n-k}^{+}=(n+2)\sum_{k=0}^{n}B_k^{+}B_{n-k}^{+},\] then we know that
\begin{equation}
\sum_{n=0}^{\infty}(n+2)\sum_{k=0}^{n}B_k^{+}B_{n-k}^{+}=4\zeta(2)\zeta(3)({\rm B}).
\end{equation}
Theorem 3.5 also tells us \[\sum_{k=1}^{\infty}\binom{k+2}{3}B_{k-1}^{+}=4\zeta(5)({\rm B}),\]
together with (3.48) we have \begin{equation}
\sum_{n=0}^{\infty}\sum_{k=0}^{n}\binom{n+2}{k}B_k^{+}B_{n-k}^{+}=2\zeta(2)\zeta(3)+2\zeta(5)({\rm B}).
\end{equation}
If we express this relation as sequence addition, we have \begin{align*}
&(2B_0^{+},6B_0^{+}B_{1}^{+},\cdots,(n+2)\sum_{k=0}^{n}B_k^{+}B_{n-k}^{+},\cdots)+(0,B_0^{+},\cdots,\binom{n+2}{3}B_{n-1}^{+},\cdots)\\
&=(2B_0^{+},4B_0^{+},\cdots,2\sum_{k=0}^{n}\binom{n+2}{k}B_k^{+}B_{n-k}^{+},\cdots).
\end{align*}
The reader can check this is  right. In fact we know
\begin{align}&\sum_{k=0}^{\infty}\frac{B_k^{+}}{k!}z^{k}=\frac{-z}{e^{-z}-1},\\
&\sum_{k=0}^{\infty}\frac{B_k^{+}}{(k+2)!}z^{k+2}=\int_{0}^{z}\int_{0}^{t}\frac{-m}{e^{-m}-1}{\rm d}m{\rm d}t.\end{align}
Multiplying above two series, we get \begin{align}
\sum_{n=0}^{\infty}\sum_{k=0}^{n}z^{2}\frac{B_k^{+}B_{n+2-k}^{+}}{k!(n+2-k)!}z^{n}=\sum_{n=0}^{\infty}\sum_{k=0}^{n}\binom{n+2}{k}B_k^{+}B_{n+2-k}^{+}\frac{z^{n+2}}{(n+2)!}.
\end{align}
We want to know the Borel sum of $\sum_{n=0}^{\infty}\sum_{k=0}^{n}\binom{n+2}{k}B_k^{+}B_{n+2-k}^{+}$, which is	 equvialent to \begin{equation}
\int_{0}^{\infty}e^{-t}\frac{-t}{e^{-t}-1}(\int_{0}^{t}\int_{0}^{x}\frac{-m}{e^{-m}-1}{\rm d}m{\rm d}x){\rm d}t.
\end{equation}
First we need to figure out the double integral. Applying series expansion, the first integral is 
 \begin{align*}\int_{0}^{x}\frac{-m}{e^{-m}-1}{\rm d}x&=\int_{0}^{x}(m+m\sum_{k=1}^{\infty}e^{-km}){\rm d}x \\
 &=\frac{1}{2}x^{2}-\sum_{k=1}^{\infty}(\frac{xe^{-kx}}{k}+\frac{e^{-kx}}{k^{2}}-\frac{1}{k^{2}})\\
&=\frac{1}{2}x^{2}+\zeta(2)-\sum_{k=1}^{\infty}(\frac{xe^{-kx}}{k}+\frac{e^{-kx}}{k^{2}}).\end{align*}
The double integral is 
\begin{align*}\int_{0}^{t}\int_{0}^{x}\frac{-m}{e^{-m}-1}{\rm d}m{\rm d}x&=\int_{0}^{t}\frac{1}{2}x^{2}{\rm d}x+\int_{0}^{t}\zeta(2){\rm d}x-\int_{0}^{t}\sum_{k=1}^{\infty}(\frac{xe^{-kx}}{k}+\frac{e^{-kx}}{k^{2}}){\rm d}x \\
	&=\frac{1}{6}t^{3}+\zeta(2)t+\sum_{k=1}^{\infty}(\frac{te^{-kt}}{k^{2}}+\frac{2e^{-kt}}{k^{3}})-2\zeta(3).
\end{align*} 
From these we know \begin{align*}
&\int_{0}^{\infty}e^{-t}(\frac{-t}{e^{-t}-1}\int_{0}^{t}\int_{0}^{x}\frac{-m}{e^{-m}-1}{\rm d}m{\rm d}x){\rm d}t\\
&=-\frac{1}{6}\int_{0}^{\infty}\frac{e^{-t}t^{4}}{e^{-t}-1}{\rm d}t-\zeta(2)\int_{0}^{\infty}\frac{e^{-t}t^{2}}{e^{-t}-1}{\rm d}t+2\zeta(3)\int_{0}^{\infty}\frac{e^{-t}t}{e^{-t}-1}{\rm d}t\\ 
&+\int_{0}^{\infty}e^{-t}(\frac{-t}{e^{-t}-1}\sum_{k=1}^{\infty}\frac{te^{-kt}}{k^{2}}){\rm d}t+\int_{0}^{\infty}e^{-t}(\frac{-t}{e^{-t}-1}\sum_{k=1}^{\infty}\frac{2e^{-kt}}{k^{3}}){\rm d}t\\
&=4\zeta(5)-2\zeta(2)\zeta(3)+2\zeta(3)\zeta(2)\\
&+\int_{0}^{\infty}e^{-t}\sum_{n=1}^{\infty}\sum_{k=1}^{n}\frac{t^{2}e^{-nt}}{k^{2}}{\rm d}t+\int_{0}^{\infty}e^{-t}\sum_{n=1}^{\infty}\sum_{k=1}^{n}\frac{2te^{-nt}}{k^{3}}{\rm d}t \\
&=4\zeta(5)+\sum_{n=1}^{\infty}\frac{2}{(n+1)^{3}}\sum_{k=1}^{n}\frac{1}{k^{2}}+\sum_{n=1}^{\infty}\frac{1}{(n+1)^{2}}\sum_{k=1}^{n}\frac{2}{k^{3}}\\
&=4\zeta(5)+2\zeta(2)\zeta(3)-2\zeta(5)\\
&=2\zeta(5)+2\zeta(2)\zeta(3).
\end{align*}

Miki \cite{miki1978relation} proved (3.36) in 1978, and it also has its divergent explanation. Let \begin{align}
	\Lambda(z)=\sum_{k=0}^{\infty}\beta_{k+1}^{+}z^{k}, \end{align}where $\beta_{k}^{+}=\beta_{k},\ k\neq1,$ and $ \beta_1^{+}=-\beta_1,$ 
then the Borel transform of $\Lambda(z)$ is 
\begin{equation}
B\Lambda(z)=\sum_{k=0}^{\infty}\frac{\beta_{k+1}^{+}}{k!}z^{k}=\frac{-1}{e^{-z}-1}-\frac{1}{z}.
\end{equation}
The Borel sum of it at $z=1$ is 
\begin{align}\int_{0}^{\infty}e^{-t}(\frac{-1}{e^{-t}-1}-\frac{e^{-t}}{t}){\rm d}t=\int_{0}^{\infty}(\frac{1}{e^{t}-1}-\frac{1}{e^{t}t}){\rm d}t=\gamma,\end{align}
where $\gamma$ is Euler-Mascheroni constant. Hence we have \begin{equation}
\sum_{k=1}^{\infty}\beta_{k}^{+}=\gamma({\rm B}).
\end{equation}
The relation between Miki's identity and divergent series will be an exploring exercise.

For more relations related to Bernoulli numbers, see \cite{agoh2014miki}.
\section{values of riemann zeta function}
 Euler found beautiful formulas for $\zeta(2k)$, which is 
 \begin{equation}
 \zeta(2n)=\frac{(-1)^{n-1}B_{2n}(2\pi)^{2n}}{2(2n)!}=\frac{(-1)^{n-1}2^{2n-1}B_{2n}\pi^{2n}}{(2n)!}.
 \end{equation}
 But the odd values of $\zeta$ are more mysterious than we think. In 1978 Ap\'ery \cite{apery1979irrationalite}  proved the irrationality of $\zeta(3)$, and it was later shown in \cite{rivoal2000fonction} that infinitely many of the odd values must be irrational. We hope there are explicit formulas as in the case of $\zeta(2n)$. However, it's untraceable. From our divergent series, we can see the values of $\zeta$ seem to have the same level. For example,
 if we write Theorem 3.5 as vector dot product, i.e.
 \begin{align*}
 &(1,-1,1,-1,\cdots)\cdot(B_0,B_1,B_2,B_3,\cdots)=\Gamma(2)\zeta(2) \\
 &(1,-2,3,-4,\cdots)\cdot(B_0,B_1,B_2,B_3,\cdots)=\Gamma(3)\zeta(3)\\
 &(2,-6,12,-20,\cdots)\cdot(B_0,B_1,B_2,B_3,\cdots)=\Gamma(4)\zeta(4) \\
 &(6,-24,60,-120,\cdots)\cdot(B_0,B_1,B_2,B_3,\cdots)=\Gamma(5)\zeta(5)\\
&(24,-120,360,-840, \cdots)\cdot(B_0,B_1,B_2,B_3,\cdots)=\Gamma(6)\zeta(6)\\
&\vdots
\end{align*}
and we also have
\begin{align*}
	&(1,-1,1,-1,\cdots)=\frac{1}{2} ({\rm B})\\
	&(1,-2,3,-4,\cdots)=\Gamma(2)\frac{1}{2^{2}} ({\rm B})\\
	&(2,-6,12,-20,\cdots)= \Gamma(3)\frac{1}{2^{3}} ({\rm B})\\
	&(6,-24,60,-120,\cdots)=\Gamma(4)\frac{1}{2^{4}} ({\rm B})\\
	&(24,-120,360,-840,\cdots)=\Gamma(5)\frac{1}{2^{5}} ({\rm B})\\
	& \vdots
\end{align*}
 It seems that there is a pattern in these values. If we look at Euler's formulas, it's very likely that  $\zeta(2n+1)$ has connection with $\pi^{2n+1}$. For instance one may want to find the value of $3\zeta(3)/{2\pi^{3}}$, but it's still a hard work.
 
 Another attempt is using the product of divergent series. For example,
 \begin{align}
 &\sum_{k=0}^{\infty}6B_{k}^{+}=\pi^{2},\\
 &\sum_{k=0}^{\infty}(k+1)B_{k}^{+}=2\zeta(3),\\
 & \sum_{k=0}^{\infty}15(k+2)(k+1)B_k^{+}=\pi^{4}.
 \end{align}
 We can deduce a divergent series which multiplyed by (3.32) equals (4.3).
 Let $\sum_{k=0}^{\infty}b_k$ denote the required series, then \begin{align*}
b_0&=B_0^{+}, \\
 b_1&=B_1^{+}, \\
 \sum_{k=0}^{n}B_{k}^{+}b_{n-k}&=(n+1)B_{n}^{+}.
 \end{align*}
Multiply (4.2) to $\sum_{k=0}^{\infty}a_k$ and make it equal $\pi^{4}$, then we get  
\begin{align*}
a_0&=5B_0^{+},\\
a_1&=10B_1^{+},\\
\sum_{k=0}^{n}6B_{n-k}^{+}a_{n-k}&=15(n+2)(n+1)B_{n}^{+}.
	\end{align*}
 \begin{remark}
 	It seems very hard to find the values of $\zeta(2n+1)$ by using product and addition on divergent series. 
 \end{remark}

\bibliographystyle{plain}
\bibliography{refere}

\begin{thebibliography}{10}

\bibitem{agoh2014miki}
T~Agoh.
\newblock On the miki and matiyasevich identities for bernoulli numbers.
\newblock {\em Integers}, 14:A17, 2014.

\bibitem{ahlfors1953complex}
Lars~V Ahlfors.
\newblock Complex analysis: an introduction to the theory of analytic functions
  of one complex variable.
\newblock {\em New York, London}, page 177, 1953.

\bibitem{apery1979irrationalite}
Roger Ap{\'e}ry.
\newblock Irrationalit{\'e} de $\zeta$ (2) et $\zeta$ (3).
\newblock {\em Ast{\'e}risque}, 61(11-13):1, 1979.

\bibitem{apostol2013introduction}
Tom~M Apostol.
\newblock {\em Introduction to analytic number theory}.
\newblock Springer Science \& Business Media, 2013.

\bibitem{bagni2009inf}
Giorgio~T Bagni et~al.
\newblock Infinite series from history to mathematics education.
\newblock {\em Colecci{\'o}n Digital Eudoxus}, 1(5), 2009.

\bibitem{borel1975lectures}
Emile Borel and Georges Bouligand.
\newblock {\em Lectures on divergent series}.
\newblock US, Energy Research and Development Administration, 1975.

\bibitem{gessel2005miki}
Ira~M Gessel.
\newblock On miki's identity for bernoulli numbers.
\newblock {\em Journal of Number Theory}, 110(1):75--82, 2005.

\bibitem{he2018new}
Chenfeng He.
\newblock A new way to understanding $\zeta(1-k)={-B_k/k}$.
\newblock {\em arXiv preprint arXiv:1811.09226}, 2018.

\bibitem{Ragib2012evaluating}
Ragib~Zaman (https://math.stackexchange.com/users/14657/ragib zaman).
\newblock Evaluating the definite integral $\int_0^\infty
  \frac{\mathrm{e}^x}{\left(\mathrm{e}^x-1\right)^2}\,x^n \,\mathrm{d}x$.
\newblock Mathematics Stack Exchange.
\newblock URL:https://math.stackexchange.com/q/117532 (version: 2012-03-07).

\bibitem{matiyasevich1997identities}
Yu~Matiyasevich.
\newblock Identities with bernoulli numbers.
\newblock {\em Preprint (http://logic. pdmi. ras.
  ru/yumat/Journal/Bernoulli/bernulli. html)}, 1997.

\bibitem{miki1978relation}
Hiroo Miki.
\newblock A relation between bernoulli numbers.
\newblock {\em Journal of Number Theory}, 10(3):297--302, 1978.

\bibitem{M1994remark}
Ján Mináě.
\newblock A remark on the values of the riemann zeta function.
\newblock {\em Expositiones Mathematicae}, 12, 01 1994.

\bibitem{murty2000simple}
M~Ram Murty and Marilyn Reece.
\newblock A simple derivation of {$\zeta(1-K)=-B_K/K$}.
\newblock {\em Funct. Approx. Comment. Math}, 28:141--154, 2000.

\bibitem{rivoal2000fonction}
Tanguy Rivoal.
\newblock La fonction z{\^e}ta de riemann prend une infinit{\'e} de valeurs
  irrationnelles aux entiers impairs.
\newblock {\em Comptes Rendus de l'Acad{\'e}mie des Sciences-Series
  I-Mathematics}, 331(4):267--270, 2000.

\end{thebibliography}

\end{document}